\documentclass[12pt,oneside,a4paper,leqno,final]{amsart}
\usepackage[all,cmtip]{xy}
\CompileMatrices

\newdir{ >}{!/10pt/@{ }*@{>}}
\newdir{< }{!/10pt/@{ }*@{<}}

\usepackage{graphicx}
\usepackage{amsfonts,amssymb}
\usepackage[centertags]{amsmath}
\usepackage{amsthm} %[2000/10/26]
\newcounter{theorems}

\theoremstyle{plain}

\swapnumbers

\numberwithin{equation}{section}

\newtheoremstyle{par}% name
     {\topsep}%      Space above (it was 0pt!!!)
     {\topsep}%      Space below (it was 0pt!!!)
     {\itshape}%         Body font
     {}%         Indent amount (empty = no indent, \parindent = para indent)
     {\bfseries}% Thm head font
     {}%        Punctuation after thm head
     {.5em}%     Space after thm head: " " = normal interword space;
     {}%         Thm head spec (can be left empty, meaning `normal')

\newtheoremstyle{parrm}% name
     {\topsep}%      Space above (it was 0pt!!!)
     {\topsep}%      Space below (it was 0pt!!!)
     {\normalfont}%         Body font
     {}%         Indent amount (empty = no indent, \parindent = para indent)
     {\itshape}% Thm head font
     {}%        Punctuation after thm head
     {.5em}%     Space after thm head: " " = normal interword space;
     {}%         Thm head spec (can be left empty, meaning `normal')

\newtheorem{theo}{\sc Theorem}[section]      % numbered within each section
\newtheorem{nr}[theo]{\sc Lemma}            % numbered along with Theorem
\newtheorem{propo}[theo]{\sc Proposition}     % numbered along with Theorem
\newtheorem{remark}[theo]{\sc Remark}           % numbered along with Theorem
\makeatletter
\def\tagform@#1{\maketag@@@{\ignorespaces#1\unskip\@@italiccorr}}

\makeatother

\usepackage{paralist}
\usepackage{comment}
\usepackage{subfigure}
\newcommand{\RR}{\mathbb{R}}
\newcommand{\CC}{\mathbb{C}}
\newcommand{\ZZ}{\mathbb{Z}}

\newcommand{\from}{\colon}

\newcommand{\R}{\mathbb R}

\newcommand{\n}[1]{{\bf #1}}

\usepackage{bm}
\newcommand{\simbolovettore}[1]{{\boldsymbol{#1}}}

\newcommand{\vp}{\simbolovettore{p}}
\newcommand{\vq}{\simbolovettore{q}}

\newcommand{\vs}{\simbolovettore{s}}

\newcommand{\vw}{\simbolovettore{w}}

\newcommand{\vz}{\simbolovettore{z}}
\newcommand{\zero}{\boldsymbol{0}}

\newcommand{\norm}[1]{\protect\left\protect\Vert\protect#1\protect\right\protect\Vert}

\newcommand{\parabolicmanifold}{\ensuremath{P}}

\newcommand{\perron}{\widehat P}

\usepackage{multirow,rotating}

\begin{document}
\pagenumbering{arabic}

\title{%
On the dihedral $n$-body problem }

\author{Davide L.~Ferrario and Alessandro Portaluri}

\date{%
\today}

\begin{abstract}
Consider $n=2l\geq 4$ point particles with equal masses in space,
subject to the following symmetry constraint: at each instant they
form an orbit of the dihedral group $D_l$, where $D_l$ is the group
of order $2l$ generated by two rotations of angle $\pi$ around two
secant lines in space meeting at an angle of $\pi/l$. By adding a
homogeneous potential of degree $-\alpha$ for $\alpha \in (0,2)$
(which recovers the gravitational Newtonian potential), one finds a
special $n$-body problem with three degrees of freedom, which is a
kind of generalisation of Devaney isosceles problem, in which all
orbits have zero angular momentum. In the paper we find all the
central configurations and we compute the dimension of the
stable/unstable manifolds.

\noindent {\em MSC Subject Class\/}: Primary 70F10; Secondary 37C80.
\vspace{0.5truecm}

\noindent {\em Keywords\/}: Dihedral $n$-body problem, McGehee
coordinates, central configurations.
\end{abstract}

\maketitle

%%%=========================================================================
\section{Introduction}
\label{sec:intro}

The goal of this paper is to compute all the central configurations
and the dimension of the stable/unstable manifolds for the dihedral
symmetric $n$-body problem in space under the action of a
homogeneous potential of degree $-\alpha$. For the Newtonian
potential this problem is a kind of generalisation of Devaney planar
isosceles three body problem \cite{Devaney80,Devaney81}.%%, following
%%Moeckel's approach to the study of the three body problem in space
%%\cite{Moeckel81,Moeckel83}.
The dihedral problem is a special case
of the full $n$-body problem which reduces to a Hamiltonian system
with three degrees of freedom. Briefly, one takes $n=2l\geq 4$ equal
masses whose initial position and velocity are symmetric with
respect to the dihedral group of rotations $D_l\subset SO(3)$. So
the masses form a (possibly degenerate and non-regular) antiprism in
space (and they are vertices of  two symmetric parallel $l$-gons).
Because of the symmetry of the problem, the masses will remain in
such a configuration for all time. Hence we have a system with only
three degrees of freedom. For $l=2$, the four bodies are at vertices
of a tetrahedron, and the problem has been studied in a series of
papers by Delgado and Vidal \cite{Vidal99,DelgadoVidal99}. The main
tool is the use of McGehee coordinates introduced in
\cite{McGehee74} but
for a general homogeneous potential of degree $-\alpha$. %%and with a
%%slight change: we consider McGehee coordinates not only for studying
%%the behaviour of solutions passing close to a total collision, but
%%also for parabolic orbits connecting central configurations,
%%projecting the full phase space to a codimension $1$ subspace.
We replace the singularity due to total collapse with an invariant
immersed manifold in the full phase space usually called \emph{total
collision manifold} which is the immersion of the parabolic manifold
of the projected phase space. We explicitly compute all central
configurations for this problem and show that just three types can
arise: a planar regular $2l$-gon, a regular $l$-gonal prism and a
$l$-gonal anti-prism. %%In particular, we prove the existence of some
%%heteroclinic connections between these central configurations and
%%some asymptotic sets on the projection of the regularised flow on
%%the parabolic manifold. Moreover by using this flow we are able to
%%establish some global results on the behaviour of solutions. We
%%discuss the qualitative behaviour of orbits which reach or come
%%close to the total collision and of those orbits which start from
%%total collapse.

The motivation in order to study this kind of problem is twofold.
From one side this problem is difficult enough to put on evidence
some chaotic behaviour of the full $n$-body problem and at the same
time it is simple enough to carry out some explicit computations.
From the other side the interest in this kind of problem is due to
the fact that it includes a lot of other problems with two or three
degrees of freedom studied in the past decades. The literature is
quite broad and we limit ourself to quote only some closest results;
among the others is the tetrahedral four body problem without and
with rotation, studied respectively in \cite{DelgadoVidal99} and
\cite{Vidal99}, the rectangular four body problem studied by Sim\'o
and Lacomba in \cite{SimoLacomba82}. %%What makes this problem harder
%%than  the tetrahedral four body problem is, also, the need of
%%computing all  central configurations.

\subsection*{Acknowledgements} We are very grateful to the anonymous
referees for their suggestions, comments and criticism which greatly
improved the manuscript.

%%%=========================================================================
\section{McGehee coordinates, projections and regularisations}
\label{sec:mcgehhecoord}

Let $V = \R^d$ denote the Euclidean space of dimension $d$ and $n
\geq 2$ an integer. Let $0$ denote the origin $0\in\R^d$. Let $m_1,
\dots, m_n$ be $n$ positive numbers (which can be thought as
masses). The configuration space of $n$ point particles with masses
$m_i$ respectively and center of mass in $0$ can be identified with
the subspace of $V^n$ consisting of all points $\vq=(\vq_1, \dots
\vq_n)\in V^n$ such that $\sum_{i=1}^n m_i\vq_i=0$. Let $\n{n}$
denote the set $\{1, \dots, n\}$ of the first $n$ positive integers.
For each pair of indexes $i,j \in \n{n}$ let $\Delta_{i,j}$ denote
the collision set of the $i$-th and $j$-th particles $\Delta_{i,j} =
\{\vq \in X|\vq_i=\vq_j\}$. Let $\Delta= \cup_{i,j} \Delta_{i,j}$ be
the {\em collision set\/}. %%The set of collision-free configurations
%%is denoted by $\hat X = X \backslash \Delta$.

Let $X \subset V^n$ be an open cone ($\R X =X$) and let $\alpha>0$
be a given positive real number. We consider the potential function
(the opposite of the potential energy) defined by
\[
U(\vq):= \sum_{i< j} \dfrac{m_im_j}{|\vq_i - \vq_j|^{\alpha}}.
\]
%%
%%
%%Consider $n$ point particles with masses $m_i$  and positions $\vq_i
%%\in X$ for $i \in\n{n}$. (we assume $X$ to be open in $\RR^k$).
%%%Consider a point $\vq$ in the configuration space $X \subset \RR^k$
%%If $M$ is the diagonal matrix (with positive entries, which can be
%%thought as masses) and $U(\vq)$ is a positive potential,
If $M$ is the diagonal matrix, then Newton equations
\[
M \ddot \vq =  \dfrac{\partial U}{\partial\vq}
\]
can be written in Hamiltonian form as
\begin{equation}\label{eq:hamilton}
\left\{
\begin{aligned}
M \dot \vq &=  \vp  \\
\dot \vp &= \dfrac{\partial U}{\partial \vq},
\end{aligned}\right.
\end{equation}
where the Hamiltonian is $H =H(\vq,\vp) = \langle\dfrac{1}{2} M^{-1}
\vp,\vp \rangle - U(\vq)$. Then equations \ref{eq:hamilton} can be
written in polar coordinates by setting the mass norm in $V^n$
defined for every $\vq \in X$ as
\[
\norm{\vq}^2 =  \langle M \vq, \vq \rangle
\]
and suitably rescaling the momentum as follows
\[\begin{aligned}
\rho & =  \norm{\vq} \\
\vs &= \dfrac{\vq}{\rho}\\
\vz & = \rho^{\beta} \vp \qquad\qquad \text{with $\alpha = 2
\beta$.}
\end{aligned}\]
In these coordinates equations~\ref{eq:hamilton} can be read as
\begin{equation}\label{eq:mcghee1}
\left\{
\begin{aligned}
\rho' &  =  \langle \vz , \vs\rangle \rho \\
\vs' &=   M^{-1}\vz - \langle \vz,\vs\rangle  \vs  \\
\vz' &= \beta \langle\vz, \vs\rangle \vz + \dfrac{\partial
U}{\partial \vq}(\vs),
\end{aligned}
\right.
\end{equation}
where the time has been rescaled by $dt = \rho^{1+\beta} d\tau$
(that is, $\dfrac{d}{d\tau} = \rho^{1+\beta}\dfrac{d}{dt}$); now the
energy can be written as
\begin{equation}\label{eq:energy}
H = \dfrac{1}{2} \rho^{-\alpha} \langle M^{-1}\vz,\vz\rangle  -
\rho^{-\alpha} U(\vs) = \rho^{-\alpha} \left( \dfrac{1}{2}  \langle
M^{-1}\vz,\vz\rangle  - U(\vs) \right).
\end{equation}
Let $k:=dn$ and let us consider the projection $(\vq,\vp) \mapsto
(\vs,\vz)$ from the full phase space $X\times  \RR^k$ to the reduced
space $S^{k-1}\times \RR^k$ (which is the trivial $\RR^k$-bundle  on
the ellipsoid $S^{k-1}$)
\[
X  \times \RR^k \to S^{k-1} \times \RR^k.
\]
In McGehee coordinates it is easy to see that the flow on $X \times
\RR^k$ can be projected to $S^{k-1} \times \RR^k$, that is
\begin{equation}\label{eq:mcghee2}
\left\{
\begin{aligned}
\vs' &=   M^{-1}\vz - \langle \vz,\vs\rangle  \vs  \\
\vz' &= \beta  \langle\vz, \vs\rangle \vz + \dfrac{\partial
U}{\partial \vq}(\vs).
\end{aligned}
\right.
\end{equation}
Also, being $X$ a cone, it is a cone on its $(k-1)$-dimensional
intersection with the ellipsoid $S^{k-1}$, which we will denote
simply by $S= S^{k-1} \cap X$. We define the \emph{parabolic
manifold} as the projection of all zero-energy orbits (or,
equivalently, of the zero-energy submanifold of $X\times \RR^k$) in
$S\times \RR^k$, that is
\[
\parabolicmanifold :=
\{ (\vs,\vz) \in S\times \RR^k : \frac{1}{2}  \langle
M^{-1}\vz,\vz\rangle  = U(s) \} \subset S^{k-1} \times \RR^k.
\]
Its dimension is $\dim S + k -1 = 2k - 2$. This is also the
projection of McGehee total collision manifold (see
\cite{McGehee74,Devaney80,Moeckel81,Moeckel83}); the manifold of
$(\vs,\vz)$ here is not considered as embedded in the space of
$(\rho,\vs,\vz)$ with $\rho=0$. By the form of equation
\ref{eq:mcghee1}, it is easy to prove the following proposition.
\begin{nr}
Solutions of~\ref{eq:mcghee2}  in $S \times \RR^k$ are projections
of solutions of~\ref{eq:mcghee1}. The parabolic manifold
$\parabolicmanifold$ is invariant for the flow of \ref{eq:mcghee2},
and solutions in $\parabolicmanifold$ can be lifted to $X\times
\RR^k$ by integrating the equation $\rho'/\rho   =  \langle \vz ,
\vs\rangle $.
\end{nr}
The parabolic manifold $\parabolicmanifold$ is the boundary of the
$(2k-1)$-dimensional \emph{elliptic} and \emph{hyperbolic}
manifolds, defined as
\[
\text{Elliptic} = \{ (\vs,\vz) \in S\times \RR^k : \dfrac{1}{2}
\langle M^{-1}\vz,\vz\rangle  < U(s) \} \subset S^{k-1} \times
\RR^k.
\]
\[
\text{Hyperbolic} = \{ (\vs,\vz) \in S\times \RR^k : \dfrac{1}{2}
\langle M^{-1}\vz,\vz\rangle  > U(s) \} \subset S^{k-1} \times
\RR^k.
\]
They are again invariant (even if the function $\dfrac{1}{2}
\langle M^{-1}\vz,\vz\rangle$ is not an invariant of the flow in
$S^{k-1}\times \RR^k$), and correspond to projection of
elliptic/hyperbolic orbits (that is, orbits with negative/positive
energy).  In fact, any fixed-energy (negative/positive) surface is
homeomorphic to the elliptic/hyperbolic manifold. Given a solution
of \ref{eq:mcghee2} in the elliptic or hyperbolic manifolds, for
each energy value $h$ the lifted solutions in $X\times \RR^k$ can be
found simply by applying \ref{eq:energy} as
\begin{equation}\label{eq:energyrelation}
\rho^\alpha= \dfrac{\langle M^{-1}\vz,\vz\rangle  - 2U(s)}{2h}.
\end{equation}
The  parabolic manifold $\parabolicmanifold$ is fiberwise
homeomorphic to a trivial $(k-1)$-sphere bundle on $S\subset
S^{k-1}$.

The next change of coordinates, due to McGehee \cite{McGehee74}
(with a reference to Sundman \cite{sundman}), is needed for defining
the Sundman--Lyapunov coordinate $v$ and for the regularisation of
the parabolic manifold $\parabolicmanifold$. Let $v,\vw\in \RR\times
\RR^k$ be defined by
\[
\begin{cases}
v& =\langle \vz, \vs \rangle \\
\vw & =M^{-1} \vz - \langle \vz, \vs \rangle \vs.
\end{cases}
\]
Then $\vz = vM\vs +M\vw$ and  $\langle \vw, M\vs \rangle = 0$, and
equations \ref{eq:mcghee2} can be replaced by
\begin{equation}\label{eq:mcghee3}
\left\{
\begin{aligned}
v' &= \norm{\vw}^2 + \beta v^2 - \alpha U(\vs)\\
\vs' &=  \vw \\
\vw' &= -\norm{\vw}^2 \vs +(\beta-1) v\vw +  M^{-1}\nabla_{\vs}
U(\vs),
\end{aligned}
\right.
\end{equation}
where $\nabla_{\vs}$ denotes covariant derivative, i.e.~the
component of the gradient tangent to the inertia ellipsoid
$\norm{\vq}=1$:
\[
\nabla_\vs U =
% \dfrac{\partial U}{\partial \vq}(\vs)   -
%  \langle \dfrac{\partial U}{\partial\vq}(\vs),  \vs \rangle \vs
% =
\dfrac{\partial U}{\partial \vq}(\vs)    + \alpha U(\vs) M \vs.
\]
The parabolic manifold $\parabolicmanifold$ is then defined by the
equation
\[
v^2 + \norm{\vw}^2   =2 U(\vs).
\]
The trivial bundle $S\times \RR^k$ is simply decomposed as the sum
of the normal bundle $(\vs,v)$ of $S$ in $\RR^k$ and the tangent
bundle $TS$ (with coordinates $(\vs,\vw)$). By the first equation in
\ref{eq:mcghee3}
\[
v'  = \norm{\vw}^2 + \beta v^2 - \alpha U(\vs)
 = (1-\beta) \norm{\vw}^2 +
 \alpha \left( \dfrac{1}{2} ( \norm{\vw}^2 + v^2 ) - U(\vs) \right),
\]
can be deduced the well-known fact that for $0<\alpha <2$,  $v$ is a
Lyapunov function on the flow in the parabolic and hyperbolic
manifolds, and therefore the flow is dissipative (gradient-like).
Moreover, the \emph{equilibrium points} in \ref{eq:mcghee3} are the
projections of the equilibrium points of \ref{eq:mcghee1} (and the
projection is one-to-one in the parabolic manifold), which can be
found as solutions of
\begin{equation}\label{eq:centralconf}
\left\{
\begin{aligned}
v^2 & =  2 U(\vs)\\
\nabla_{\vs} U(\vs) & =   \zero \\
\vw  &=   \zero.
\end{aligned}
\right.
\end{equation}
Hence all equilibrium points belong to the parabolic manifold
$\parabolicmanifold$. The constant solution in a central
configuration $\bar \vs$ with $v^2 = 2U(\bar\vs)$ can be lifted to
the full space as a homotetic parabolic orbit by  integrating (back
to the real time coordinate)
\[
\dot \rho = \pm \rho^{-\beta}  \sqrt{ 2U(\bar \vs)} \implies
\rho (t) = \left( \pm(1+\beta) \sqrt{2U(\bar \vs)} t  \right) ^{
1/\left( 1+\beta \right) },
\]
assuming the total collision occurs at $t=0$ (the $+$ sign yields an
ejection solution, the $-$ sign yields a collision solution). More
generally, \emph{homotetic solutions} (i.e. $\vs'=0$, with $\vs(t)
\equiv \bar\vs$) can be found in the hyperbolic and elliptic
manifolds by setting in equations \ref{eq:mcghee3} $\vs' = \vw = 0$,
and therefore by integrating the single equation
\[
v' = \beta v^2 - \alpha U(\bar\vs)
\]
and then lifting the solution found to the full space using the
energy relation \ref{eq:energyrelation}. The graphs of homotetic
solutions are straight lines contained in the normal bundle of $S$
in $S\times \RR^k$.

\section{The dihedral $2n$-body problem}
\label{sec:equivsetup} 
\begin{figure}
\centering
\subfigure[$l=2$]{%
\includegraphics[width=0.480\textwidth]{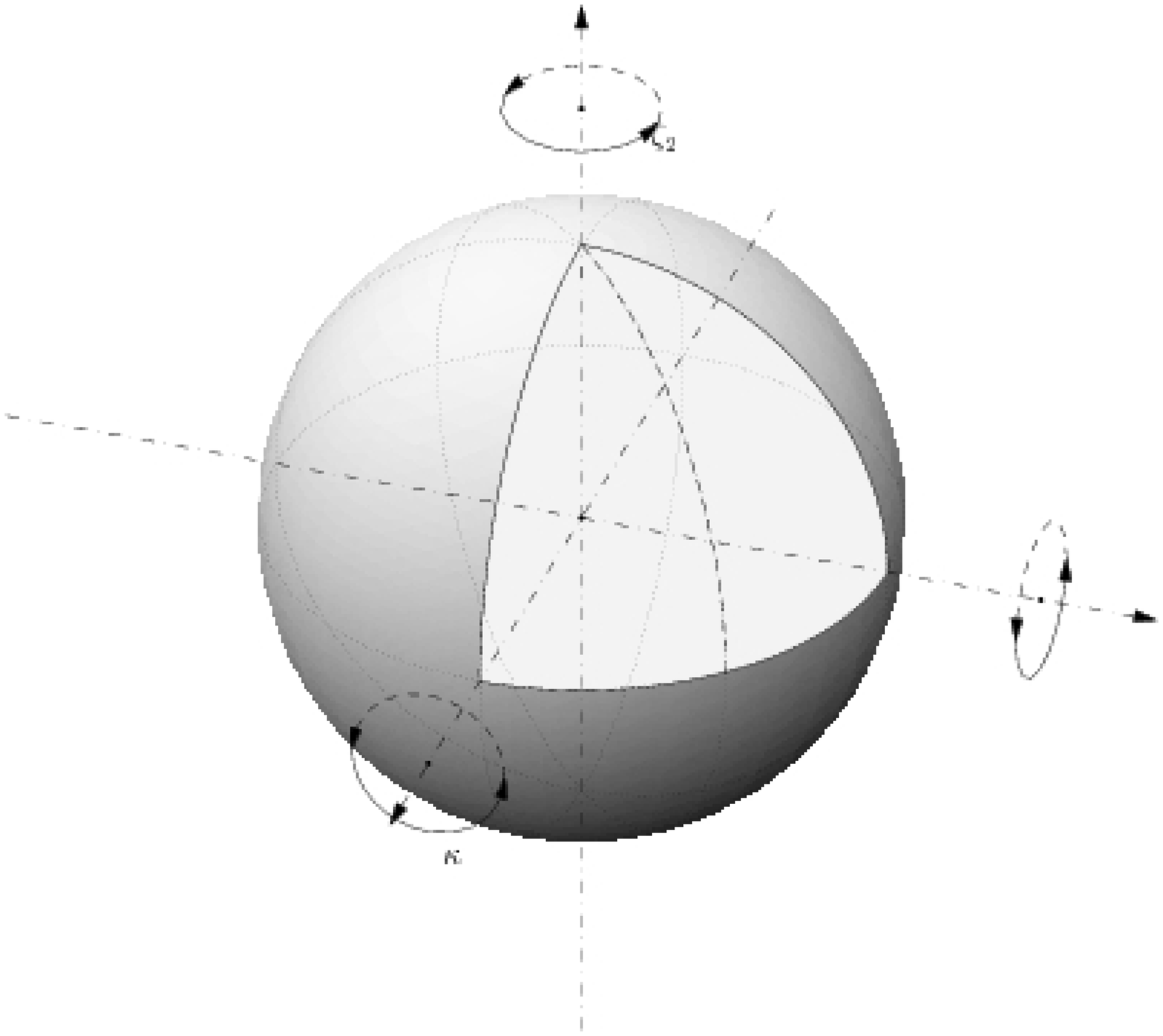}
% \caption{Fundamental domain for the dihedral action for $l=2$}
\label{fig:fund4} }
\subfigure[$l=3$]{%
% \caption{Fundamental domain for the dihedral action for $l=3$}
\includegraphics[width=0.480\textwidth]{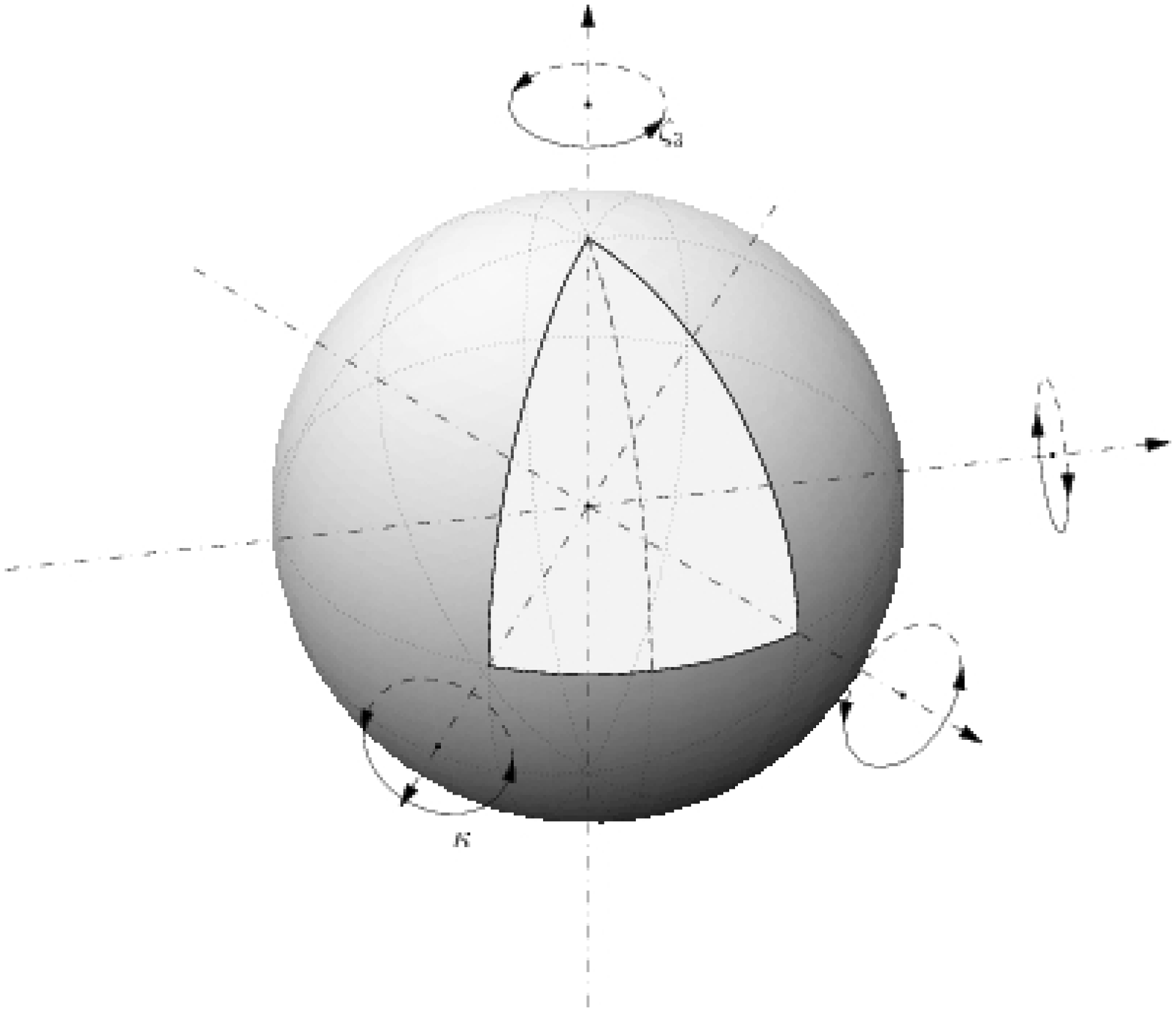}
\label{fig:fund} } \caption{Dihedral groups $D_l$, with the upper
half of the fundamental domains in white.} \label{fig:fun}
\end{figure}

Let $\RR^3 \cong \CC \times \RR$ be endowed
with coordinates $(z,y)$, $z\in \CC$, $y\in \RR$. For $l\geq 1$, let
$\zeta_l$ denote the primitive root of unity $\zeta_l = e^{2\pi
i/l}$; the \emph{dihedral} group $D_{l}\subset SO(3)$ is the group
of order $2l$ generated by the rotations
\[
\zeta_l\from (z,y) \mapsto (\zeta_l z, y ) \text{\ and \ }
\kappa\from (z,y) \mapsto (\overline{z},-y),
\]
where $\overline{z}$ is the complex conjugate of $z$. The
non-trivial elements of $D_l = \langle \zeta_l,\kappa \rangle$ are
the  $l-1$ rotations around the $l$-gonal axis $\zeta_l^j$,
$j=1,\ldots,l-1$ and  the $l$ rotations of angle $\pi$ around the
$l$ digonal  axes orthogonal to the $l$-gonal axis (see figure
\ref{fig:fun}) $\zeta_l^j\kappa$, $j=1,\ldots ,l$. In
figures~\ref{fig:fund4} and~\ref{fig:fund} one can find the
upper-halves of the fundamental domains for the action of $D_l$
restricted on the unit sphere. In fact, in figure~\ref{fig:fund4}
corresponding to the dihedral four body problem, the fundamental
domain is represented by an octant of the shape sphere while
figure~\ref{fig:fund} represent the fundamental domain on the shape
sphere for the dihedral six body problem.
%%
%%\begin{figure}\centering
%%\includegraphics[width=0.618\textwidth]{figs/fun2}
%%\caption{$G$-orbit of a generic point for the dihedral action with $l=3$
%%(6 bodies).}
%%\label{figs/fun2}
%%\end{figure}

Consider the permutation representation of $D_l$ given by left
multiplication (that is, the Cayley immersion $\sigma\from D_l \to
\Sigma_{2l}$ of $D_l$ into the symmetric group on the  $2l$ elements
of $D_l$, defined by $\sigma(g)(x) = gx$ for each $g,x\in D_l$, see
\cite{gg} for more details). The action of $D_l$ on $\RR^3$ induces
an orthogonal action on the configuration space $\RR^{6l}$ of $n=2l$
point particles $\vq_i\in \RR^3$ in the three-dimensional space. The
Newtonian potential for the $n$-body problem, homogeneous with
degree $-\alpha$ induces by restriction on the fixed subspace
$\left(\RR^{6l}\right)^{D_{l}} \cong \RR^3$ a homogeneous potential
defined for each $\vq\in \RR^3$ by
\begin{equation}\label{eq:potenziale}
U(\vq) = \sum_{g\in D_l \smallsetminus \{1\}} \left| \vq - g \vq
\right|^{-\alpha},
\end{equation}
provided we assume (without loss of generality) all masses
$m_i^{2}=1/l$. Now, the potential $U$ in \ref{eq:potenziale} can be
re-written in terms of coordinates $\vq = (z,y) \in \CC\times \RR$
as
\[\begin{aligned}
U(\vq) &= \sum_{j=1}^{l-1} | \vq - \zeta_l^j \vq |^{-\alpha} +
\sum_{j=1}^{l} | \vq - \zeta_l^j \kappa \vq |^{-\alpha} \\
&= \sum_{j=1}^{l-1} | z - \zeta_l^j z |^{-\alpha} + \sum_{j=1}^{l}
\left(  | z - \zeta_l^j \overline{z}|^2 + 4y^2 \right) ^{-\alpha/2}.
\end{aligned}
\]
By definition, for each $g\in D_{l}$, $U(g\vq)=U(\vq)$. Further
symmetries of $U$ are:
\begin{enumerate}
\item the reflection on the
plane $y=0$ (given by $h\from (z,y) \mapsto (z,-y)$),
\item
the $l$ reflections on the  planes containing the $l$-gonal axis and
one of the digonal axes,
\item and the $l$ reflections on the planes
containing the $l$-gonal axis and the points $( \zeta_l^j e^{\pi i/
l}  , 0 )$, $j = 1, \ldots, l$.
\end{enumerate}
It is not difficult to prove that these are (up to conjugacy and
multiplication with elements in $D_l$) all the elements of  the
normaliser of $D_l$ in $O(3)$. Thus we can study $U$ only in the
left-upper area of the $D_l$-fundamental domain on $S^2\subset
\RR^3$, as we have seen in figures \ref{fig:fund4} and
\ref{fig:fund}. Now, in order to simplify the expression of the
potential we introduce the variables $r$ and $\xi$ as follows. If
$y\geq 0$ and $z\neq 0$, let $r=r(z,y)$ be defined as $r = 1 +
2y^2/|z|^2 - 2y/|z|\sqrt {y^2/|z|^2 + 1}$ and $\xi =
\dfrac{\overline{z}}{z}$. Hence $r\in (0,1]$, with $r=1$ if and only
if  $y=0$, $1+r^2 = r ( 2 + 4 y^2/|z^2|)$ and therefore
\[
\begin{aligned}
| z - \zeta_l^j \overline{z}|^2 + 4y^2  &= |z|^2 \left( | 1 -
\zeta_l^j \dfrac{\overline{z}}{z} |^2 + 4 \dfrac{y^2}{|z|^2} \right)
 = \dfrac{|z|^2}{r} \left| 1 -  \zeta_l^j r \xi \right|^2.
\end{aligned}
\]
In these coordinates the potential function $U(\vq)$ can be written
as
\begin{equation}\label{eq:potentialplain}
U = |z|^{-\alpha} \left[ \sum_{j=1}^{l-1} | 1 - \zeta_l^j
|^{-\alpha} + r^{-\alpha/2} \sum_{j=1}^{l} \left| 1 - \zeta_l^j r
\xi \right| ^{-\alpha} \right].
\end{equation}
We can now state the integral representation of the potential
\ref{eq:potentialplain} proven in the Appendix \ref{sec:appendice}
(see also \cite{BanElm} and remark \ref{remark:qqq} below).
\begin{propo}
\label{propo:potentialintegral} For $\beta\in(0,1)$, $r\in (0,1]$
and $\xi\in S^1\subset \CC$ the potential $U$ can be written as
\[
U = |z|^{-\alpha} \left[ c_l  + l r^{-\beta} \dfrac{\sin( \beta
\pi)}{ \pi } \int_0^1
\dfrac{(1-t)^{-\beta} t^{\beta-1} }{ %
{(1-tr^2)}^{\beta} } \dfrac{1 - (tr)^{2l}}{|1-(tr)^l\xi^l|^2} \,dt
\right],
\]
where $c_l$ is the constant $c_l  = \sum_{j=1}^{l-1} | 1 - \zeta_l^j
|^{-\alpha}$.
\end{propo}
\proof For the proof of this result, see Appendix
\ref{sec:appendice}.\qed
\begin{remark}
The above integral representation plays a fundamental role in order
to find all the central configuration. In fact, otherwise the
expression of the potential given in formula \ref{eq:potentialplain}
is quite difficult to deal with.
\end{remark}
\subsection{Planar type central configurations}

On the unit sphere $S\subset \RR^3$ (of equation $|z|^2 + y^2 = 1$),
parametrised by $(\varphi,\theta)\in (-\pi/2,\pi/2) \times [0,2\pi)$
with $y=\sin  \varphi$ and $z = \cos\varphi e^{i\theta}$, the
(reduced to the $2$-sphere) potential reads
\begin{equation}\label{potential:sphere}
{U}(\theta,\varphi) \!\!=\!\! (2\cos\varphi)^{-\alpha} \left[
\sum_{j=1}^{l-1} \left( \sin \dfrac{j\pi}{l} \right)^{-\alpha}
\!\!\!+\!\! \sum_{j=1}^{l} \left( \sin^2 (\dfrac{j\pi}{l}
\!-\!\theta ) \!+\! \tan^2 \varphi \right) ^{-\frac{\alpha}{2}}
\right],
\end{equation}
and by Proposition \ref{propo:potentialintegral} also as
\[
{U}(\theta,r) = \left( \dfrac{1+r^2}{4r} \right)^{\beta} \left[ c_l
+ l r^{-\beta} \dfrac{\sin( \beta \pi)}{ \pi } I(r, \theta)\right]\]
\[
\textrm{where}\ \ I(r,\theta)=\int_0^1
\dfrac{(1-t)^{-\beta} t^{\beta-1} }{ %
{(1-tr^2)}^{\beta} } \dfrac{1 - (tr)^{2l}}{1+(tr)^{2l} -
2(tr)^l\cos(2 l\theta) } \,dt ,
\]
with (just for $\varphi\in [0,\pi/2)$)
\[
r = 1 + 2\tan^2\varphi  -2 \dfrac{\tan\varphi}{\cos\varphi} =
\dfrac{1-\sin\varphi}{1+\sin\varphi}
\]
and hence
\[
\sin\varphi = \dfrac{1-r}{1+r} \text{\ and \ } \cos^2\varphi =
\dfrac{4r}{1+r^2}.
\]
In spherical coordinates, the symmetry reflections  of $U$ are (up
to conjugacy)
\begin{enumerate}
\item the reflection on the horizontal
plane: $h_\varphi\from (\theta,\varphi)  \mapsto
(\theta,-\varphi)$,
\item
the reflection on the plane containing the $l$-gonal axis and the
digonal axis $h_\theta\from  (\theta,\varphi) \mapsto
(-\theta,\varphi)$
\item and the reflection on the plane
containing the $l$-gonal axis and the point $( e^{\pi i/ l}  , 0 )$,
defined as $h'_\theta\from (\theta,\varphi) \mapsto (\pi/l - \theta
,\varphi)$.
\end{enumerate}
As direct consequence of the Palais' symmetric criticality
principle, it follows that critical points of the restrictions of
the reduced potential $U$ to the $1$-spheres of such fixed planes
are critical points for the restriction of $U$ to the sphere, and
hence are central configurations for $U$. In fact, as already
observed this $1$-spheres are nothing but the spaces fixed by each
of the reflections given in $(i)$, $(ii)$ and $(iii)$. In principle
it can be exist other critical points for the restriction of the
potential $U$ to the sphere which do not lie in these fixed spaces.
However if we are able to show that out of this $1$-spheres the
derivative of the potential is bounded away from zero, we have done.

Now consider the derivative with respect to $\theta$ of $U$, which
by Proposition \ref{propo:potentialintegral} can be written as
follows
\begin{equation}\label{eq:partialtheta}
\dfrac{\partial U}{\partial \theta}  = -4 l^2 \sin(2l\theta) \dfrac{
(1+r^2)^\beta \sin (\beta \pi ) }{ \pi ( 2r )^\alpha
% r^\beta(\cos\varphi)^{\alpha}
} I(r,\theta)
\end{equation}
where $I(r,\theta)$ is  strictly positive and defined for
$(\theta,r) \neq (2k\pi/l,1)$, $k$ integer. Hence for each $r\in
(0,1]$ the derivative $\dfrac{\partial U}{\partial \theta}$ is
strictly negative for $\theta \in (0,\dfrac{\pi}{2l})$ and strictly
positive for $\theta\in (\dfrac{\pi}{2l},\dfrac{\pi}{l})$. It is
zero for $\theta= \dfrac{k\pi}{2l}$ and $r\in (0,1)$ and
$\theta=\dfrac{(2k+1)\pi}{2l}$ and $r=1$. Thus, for $\varphi=0$, we
have proved the following proposition:
\begin{nr}[Planar $2l$-gon]
\label{nr:nagon} For any $\alpha \in (0,2)$ central configurations
which are $h_\varphi$-symmetric are on the vertices $(e^{(2k+1)\pi
i/(2l) },0)$ of the regular $2l$-gon.
\end{nr}
\subsection{Prism type central configurations}
Now we have to explore the cases $\theta = k\pi /l $ and $\theta  =
(2k+1)\pi/(2l)$, which correspond respectively to prisms and
antiprisms. The derivative of \ref{potential:sphere} with respect to
$\varphi$ is
\begin{equation}\label{eq:partialvphi}\begin{aligned}
\dfrac{\partial U}{\partial \varphi} & = 2 \beta \dfrac{
\tan\varphi}{(2\cos\varphi)^\alpha} \left[ c_l - \sum_{j=1}^l
\dfrac{ \cos^2(j\pi/l - \theta) }{ \left( \sin^2(j\pi/l - \theta) +
\tan^2\varphi \right)^{\beta+1} }
\right] \\
& = 2 \beta \dfrac{ \tan\varphi}{(2\cos\varphi)^\alpha} \left[
f_\theta(\varphi) \right].
\end{aligned}
\end{equation}
The term in square brackets $f_\theta(\varphi)$ has the same sign of
$\dfrac{\partial U}{\partial \varphi}$, and since $c_l$ is a
constant and each term of the sum is strictly monotone in $\varphi$,
for each $\theta$ the function $f_\theta(\varphi)$ can vanish at
most once in the interval $(0,\pi/2)$. Since the  limit of the sum
as $\varphi\to \pi/2$ is zero and $c_l$ is positive, there will be a
unique zero in $(0,\pi/2)$ (for a fixed $\theta$) for all the values
$\theta$ such that $\lim_{\varphi \to 0} f_\theta(\varphi) < 0 $,
i.e.
\[
\lim_{\varphi \to 0} \sum_{j=1}^l \dfrac{ \cos^2(j\pi/l - \theta) }{
\left( \sin^2(j\pi/l - \theta) + \tan^2\varphi \right)^{\beta+1} }
>  c_l = \sum_{k=1}^{l-1} \left( \sin^2 (j\pi/l) \right)^{-\beta}.
\]
Now, since $ \lim_{\varphi \to 0} f_0(\varphi) = -\infty,$ there
exists a unique minimum $\hat\varphi$ for $\theta= k\pi/l$,
$k=0\ldots 2l-1$, corresponding to a prism.

\begin{nr}[Prisms]
\label{nr:prism} There are exactly $4l$ central configurations which
are $h_\theta$-symmetric (up to conjugacy), and they are precisely
on the vertices  of a prism: $(\cos\hat\varphi' e^{k\pi i/l },\pm
\sin \hat\varphi')$.
\end{nr}
We observe that in the dihedral four body problem these kind of
central configurations collapse to square type central
configurations.
\subsection{Antiprism type central configurations}
It is left to compute critical points for $\theta = (2k+1)\pi/(2l)$,
that is, to find zeroes of $f_\theta(\varphi)$ for $\theta =
\dfrac{\pi}{2l}$, or, equivalently, $h'_\theta$-symmetric central
configurations.

\begin{nr}[Antiprisms]
\label{nr:antiprism} There are exactly $2l$ central configurations
which are $h'_\theta$-symmetric (up to conjugacy) and $\varphi\neq
0$. They are on the vertices of a prism: $(\cos\hat\varphi
e^{(2k+1)\pi i/(2l) },\pm \sin \hat\varphi)$.
\end{nr}
We remark that in the four body problem the antiprism type central
configurations reduce to tetrahedral type configurations.
\begin{proof}
It suffices to show that
\[
\sum_{j=1}^l \dfrac{ \cos^2(j\pi/l - \pi/(2l)) }{ \left(
\sin^2(j\pi/l - \pi/(2l)) \right)^{\beta+1} }
>  \sum_{k=1}^{l-1} \left( \sin^2 (j\pi/l) \right)^{-\beta}.
\]
If $\lfloor l/2 \rfloor$ denotes the greatest integer $n\leq l/2$,
that is
\[
\lfloor l/2 \rfloor =
\begin{cases}
(l-1)/2 & \text{ $l$ odd} \\
l/2 & \text{$l$ even,}
\end{cases}
\]
then
\[
\sum_{j=1}^l \dfrac{ \cos^2(j\pi/l - \pi/(2l)) }{ \left(
\sin^2(j\pi/l - \pi/(2l)) \right)^{\beta+1} } = 2
\sum_{j=1}^{\lfloor l/2 \rfloor} \dfrac{ \cos^2(j\pi/l - \pi/(2l) )
}{ \left( \sin^2(j\pi/l - \pi/(2l)) \right)^{\beta+1}}.
\]
On the other hand
\[
\sum_{j=1}^{l-1} \left( \sin^2 (j\pi/l) \right)^{-\beta} = 2
\sum_{j=1}^{\lfloor l/2 \rfloor} \left( \sin^2 (j\pi/l)
\right)^{-\beta} + d_l,
\]
where
\[
d_l = \begin{cases}
1 & \text{$l$ even} \\
0 & \text{$l$ odd}
\end{cases}
\]
Now then, since
\[
\dfrac{ \cos^2 x  }{ \left( \sin^2 x \right)^{\beta+1}} = \dfrac{ 1
}{ \left( \sin^2 x \right)^{\beta+1}} - \dfrac{ 1  }{ \left( \sin^2
x \right)^{\beta}},
\]
the conclusion would follow once we could prove that
\[
 2 \sum_{j=1}^{\lfloor l/2 \rfloor} C_j  > d_l,
\]
where
\[
C_j\! = \!\dfrac{ 1  }{ \left( \sin^2 (j\pi/l - \pi/(2l))
\right)^{\beta+1}} - \dfrac{ 1  }{ \left( \sin^2 (j\pi/l - \pi/(2l))
\right)^{\beta}} -
\dfrac{1}{%
\left( \sin^2 (j\pi/l) \right)^{\beta}}.
\]
If $l=2$, it turns out that $C_1 = 2^\beta - 1$ and hence $2
\sum_{j=1}^{\lfloor l/2 \rfloor} C_j = 2C_1 > 0 = d_2$. If $l=3$,
then $C_1 = 2^\alpha ( 3 - 3^{-\beta})$, which is greater than $2$
for all $\alpha = 2\beta$, so that $2C_1 > 1 = d_1$. In general,
since $j \leq l/2$, $\sin(j\pi/l) > \sin(j\pi/l - \pi/(2l))$, and
therefore
\[\begin{aligned}
C_j & > \dfrac{ 1  }{ \left( \sin^2 (j\pi/l - \pi/(2l))
\right)^{\beta+1}} - \dfrac{ 2  }{ \left( \sin^2 (j\pi/l - \pi/(2l))
\right)^{\beta}} \\
& = \dfrac{ 1 - 2 \sin^2( j\pi/l - \pi/(2l))  }{ \left( \sin^2
(j\pi/l - \pi/(2l)) \right)^{\beta+1}}.
\end{aligned}\]
The first term is estimated by
\[\begin{split}
C_1 \geq \dfrac{ 1 - 2 \sin^2( \pi/(2l))  }{ \left( \sin^2
(\pi/(2l)) \right)^{\beta+1}}
> \dfrac{1 -  \dfrac{\pi^2}{2l^2} }{%
\left( \dfrac{\pi}{2l} \right)^{\alpha+2} } = (1 - \dfrac{\pi^2}{2
l^2}) \dfrac{2^{\alpha+2} l^{\alpha+2} }{\pi^{\alpha+2}}
\\
\geq
% (1 - \dfrac{\pi^2}{2 l^2}) \dfrac{4 l^{2} }{\pi^{2}}
\dfrac{4l^2}{\pi^2} - 2 > \dfrac{l^2}{4} - 2,
\end{split}\]
and all other terms $C_j$ with $j\geq 2$ are in any case greater
than $-1$; thus for $l\geq 4$
\[
 \sum_{j=1}^{\lfloor l/2 \rfloor} C_j
>
 C_1 - ( \lfloor \dfrac{l}{2} \rfloor - 1 )
 \geq C_1 - \dfrac{l}{2} + 1
 \geq
 \dfrac{l^2}{4} -  \dfrac{l}{2} -1 \geq  1,
\]
and thus for all $l\geq 4$ we have $2  \sum_{j=1}^{\lfloor l/2
\rfloor} C_j \geq 2  > d_l,$ which concludes the proof.
\end{proof}

\begin{figure}\centering
\includegraphics[width=0.718\textwidth]{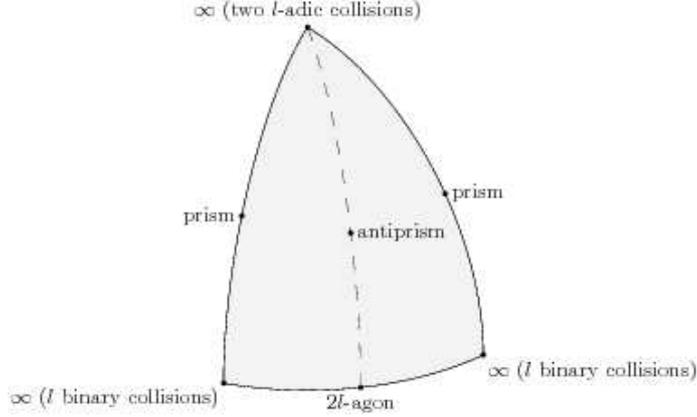}
\caption{Central configurations in the upper-half fundamental
domain.} \label{fig:centralconfigurations}
\end{figure}

Since there are no other central configurations, by
\ref{eq:partialtheta}, we can summarise the results in the following
proposition.
\begin{propo}
\label{propo:main} All central configurations in the dihedral
$2n$-body problem are symmetric for one of the three types of
reflections $h_\varphi$ Lemma \ref{nr:nagon}, $h_\theta$ Lemma
\ref{nr:prism} or $h'_\theta$ Lemma \ref{nr:antiprism}. They are
represented in the (upper-half) fundamental domain on the sphere in
figure \ref{fig:centralconfigurations}.
\end{propo}
In fact, in figure \ref{fig:centralconfigurations} it is drawn a
geodesic triangle which represents the fundamental domain on the
shape sphere for the dihedral $n$-body problem. In this figure are
shown the three types of central configurations arising in the
problem we are dealing with in the exact location together with.
Moreover we observe that due to the symmetry constraint only two
types of collisions can occur. We denoted by the name $l$-adic
collision and binary collision, meaning that in the first case two
clusters of $l$-bodies simultaneously collide, while in the second
case $l$ clusters of 2 bodies simultaneously collide. This two types
of collisions are all located on the same plane containing the
planar central configurations while the $l$-adic central
configurations can be represented in the north and south pole of the
shape sphere.

Now consider equations \ref{eq:mcghee3} in coordinates
$(\theta,\varphi)$ on the sphere: we set $w_1$ and $w_2$ such that
$\vw = w_1 \dfrac{\partial \vs}{\partial \theta} + w_2
\dfrac{\partial \vs}{\partial \varphi}$, i.e.  (since $\vs =
(\cos\varphi e^{i\theta}, \sin \varphi )$),
\[
\vw = w_1 ( i \cos\varphi e^{i\theta}, 0) + w_2 (-\sin\varphi
e^{i\theta}, \cos\varphi).
\]
Then
\[
\norm{\vw} = w_1^2\cos^2\varphi + w_2^2
\]
and
\[
\nabla_{\vs} U(\vs)  = \dfrac{1}{\cos^2\varphi} \dfrac{\partial
U}{\partial \theta} \dfrac{\partial \vs}{\partial \theta} +
\dfrac{\partial U}{\partial \varphi} \dfrac{\partial \vs}{\partial
\varphi}.
\]
Also, equations \ref{eq:mcghee3} become
\begin{equation}\label{eq:mcgheesphere}
\left\{
\begin{aligned}
v' &= w_1^2 \cos^2\varphi + w_2^2 +\beta v^2 - \alpha U(\theta,\varphi) \\
\theta' &= w_1 \\
\varphi' &= w_2 \\
w_1' &= (\beta-1)vw_1 + 2 \tan \varphi~  w_1 w_2 +
\dfrac{1}{\cos^2\varphi} \dfrac{\partial U}{\partial \theta}
\\
w_2' &= (\beta-1)v w_2 - \dfrac{1}{2} w_1^2 \sin 2\varphi  +
\dfrac{\partial U}{\partial \varphi}.
\end{aligned}
\right.
\end{equation}
The linearization at equilibrium points (central configurations)
\ref{eq:centralconf} is represented by the $5\times 5$ matrix  $L$
\[
L =
\begin{bmatrix}
2\beta v & 0 & 0 & 0 & 0 \\
0 & 0 & 0 & 1 & 0 \\
0 & 0 & 0 & 0 & 1 \\
0 & \dfrac{1}{\cos^2\varphi} \dfrac{\partial^2 U}{\partial \theta^2}
&
\dfrac{\partial^2 U}{\partial \theta \partial \varphi } & (\beta-1)v & 0 \\
0 & \dfrac{\partial^2 U}{\partial \varphi\partial\theta} &
\dfrac{\partial^2 U}{\partial\varphi^2} & 0 & (\beta-1)v \\
\end{bmatrix}.
\]
Thus the eigenvalues of the linearization can be computed in terms
of the eigenvalues of the Hessian $D^2U(\overline\vs)$  of
$U(\theta,\varphi)$:
\begin{propo}\label{propo:borat}
The eigenvalues of $L$, at a central configuration $\overline\vs$
(i.e. at the point $(\overline v, \overline \vs, \zero)$, where
$\overline v = \pm \sqrt{2U(\overline\vs)}$), are equal to the roots
$\lambda$ of the equation
\[
\lambda^2 + (1-\beta)\overline v \lambda = \gamma
\]
for each $\gamma$ eigenvalue of the Hessian $D^2U(\overline\vs)$.
\end{propo}

By elementary calculations it follows from Proposition
\ref{propo:borat} that the equilibrium points
$(\pm\sqrt{2U(\overline\vs)},\overline\vs,\zero)$ are hyperbolic
when the Hessian $D^2U$ is non-singular at $\overline{\vs}$, and
that for each positive eigenvalue $\gamma>0$ of $D^2U$ there is a
pair of real eigenvalues of $L$,  $\lambda_1>0$, $\lambda_2<0$; for
each negative eigenvalue $\gamma <0$ of $D^2U$, there are two
eigenvalues $\lambda_1$, $\lambda_2$ of $L$ with negative real part
($\lambda_1$, $\lambda_2$ are real if $d = (1-\beta)^2{\overline
v}^2 + 4\gamma > 0$ and $\lambda_1=\lambda_2$ if $d=0$).

\begin{propo}
All equilibrium points of \ref{eq:mcgheesphere} are hyperbolic.
% (with non-zero
% real part).
\end{propo}
\begin{proof}
We just need to proof that the Hessian $D^2U$ is non-singular at
$\overline\vs$, if $\vs$ is a central configuration. Since each
central configuration $\vs$ lies in the line fixed by a reflection
(which is a symmetry of $U$), the matrix $D^2U$ is diagonal at
$\overline\vs$. So the result follows once we prove that
$\dfrac{\partial^2 U}{\partial \theta^2}(\overline\vs) \neq 0 \neq
\dfrac{\partial^2U}{\partial \varphi^2}(\overline\vs)$. But by
\ref{eq:partialtheta}, since $I(r,\theta)$ is strictly positive and
regular in a neighbourhood of $\overline{s}$, $\dfrac{\partial^2
U}{\partial \theta^2}(\overline\vs) \neq 0$. By
\ref{eq:partialvphi}, the same holds for
$\dfrac{\partial^2U}{\partial \varphi^2}(\overline\vs)$.
\end{proof}

\begin{propo}\label{coro:dimensions}
The dimension of the stable (unstable) manifold of $(\overline
v,\overline\vs,\zero)$ with $\overline v = \sqrt{2 U(\overline\vs)}
> 0$ is 3 (2) if $\overline\vs$ is a $2l$-gon or a prism; it is 2
(3) if $\overline\vs$ is an antiprism. The dimension of the stable
(unstable) manifold of the point $(- \overline
v,\overline\vs,\zero)$ with $\overline v = \sqrt{2 U(\overline\vs)}
> 0$ is equal to the dimension of the unstable (stable) manifold of
$(\overline v,\overline\vs,\zero)$. The intersection of the stable
(unstable) manifold of $(\overline v,\overline\vs,\zero)$ with the
parabolic manifold $\parabolicmanifold$ has codimension 0 (1) in
$\parabolicmanifold$ if $\overline v >0$. It has codimension 1 (0)
in $\parabolicmanifold$ if $\overline v<0$.
\end{propo}
\begin{proof}
These facts follow directly from the stable/unstable manifold
theorem and the above arguments on eigenvalues of $L$. The results
are summarised in table \ref{tb:tabella1}.
\end{proof}

% \begin{table}\centering
% \begin{tabular}{c|ccccc}
% & $\overline v $&  {$\dim W^u \cap \parabolicmanifold$}
% &{$\dim W^s \cap \parabolicmanifold$}& {$\dim  W^u$}&{$\dim  W^s$}\\
% \hline \\
% $2l$-agon and prism & $>0$  &  1 & 3& 2& 3 \\
%  & $<0$  & 3&1&3&2 \\
% \hline\\
% Anti-prism  & $>0$& 2 &2& 3& 2   \\
%  &$<0$ & 2 &2& 2 & 3   \\
% \hline
% \end{tabular}
% \caption{Dimensions of stable and unstable manifolds}
% \label{tb:tabella1}
% \end{table}

\begin{table}\centering
\begin{tabular}{|l|c|cc|cc|}
\hline & $\overline v $&
\begin{sideways}$\dim  W^s$\end{sideways} &
\begin{sideways}$\dim  W^u$\end{sideways} &
\begin{sideways}$\dim W^s \cap \parabolicmanifold$\end{sideways}&
\begin{sideways}$\dim W^u \cap \parabolicmanifold$\end{sideways}
\\ \hline
\multirow{2}{*}{$2l$-gon and prism} & $>0$  &
3 & 2 & 3 & 1 \\
 & $<0$
 & 2 & 3 & 1 & 3 \\
\hline \multirow{2}{*}{anti-prism}  & $>0$&
2 &3& 2& 2   \\
 &$<0$ & 3 &2& 2 & 2   \\\hline
\end{tabular}
\caption{Dimensions of stable and unstable manifolds.}
\label{tb:tabella1}
\end{table}

\appendix
\section[Tua madre]{An integral representation for $U$}\label{sec:appendice}

The aim of this section is to give a direct proof of the integral
representation for the potential $U$ used before in order to compute
all the central configurations.

For $l\geq 2$, let $\perron_l$ denote the $l$-adic Perron-Frobenius
operator, defined on complex functions $f\from S^1\subset \CC \to
\CC$ by
\[
\begin{aligned}
\forall \xi=e^{i\theta} \in S^1, \perron_l(f)(\xi) & = \dfrac{1}{l}
\sum_{y\,:\, y^l = \xi} f(y)
%%&=
%%\dfrac{1}{l}
%%\sum_{\substack{j=0\ldots (l-1)\\ y_1: y_1^l = \zeta_l}} {f(\zeta_l^j y_1)} \\
%%\\
=  \dfrac{1}{l} \sum_{j=0}^{l-1} f( e^{\frac{i(\theta + 2j\pi)}{l} }
).
\end{aligned}
\]
For each $k\in \ZZ$,
\begin{equation}\label{eq:proprperron}
\perron_l(\xi^k)(\xi)  = \dfrac{1}{l} \sum_{y\,:\, y^l = \xi} y^k =
\begin{cases}
\xi^{k/l}   &
\text{if $k\equiv 0 \mod l$} \\
0 & \text{if $k \not\equiv 0 \mod l$.}
\end{cases}
\end{equation}
In terms of the $l$-adic Perron-Frobenius operator,  the potential
\ref{eq:potentialplain} can be written as

\begin{equation}\label{eq:potentialperron}
U = |z|^{-\alpha} \left[ c_l  + r^{-\alpha/2} l \perron_l\left(
\left| 1 - r \xi \right| ^{-\alpha} \right)(\xi^l) \right],
\end{equation}
where $c_l$ is the constant $c_l  = \sum_{j=1}^{l-1} | 1 - \zeta_l^j
|^{-\alpha}$ and $\perron_l\left( \left| 1 - r \xi \right|
^{-\alpha} \right)(\xi^l)$ denotes the function $\perron_l\left(
\left| 1 - r \xi \right| ^{-\alpha} \right)$ of argument $\xi$
evaluated at $\xi^l$. In order to compute  $\perron_l\left( \left| 1
- r \xi \right| ^{-\alpha} \right)$, we expand $|1-r\xi|^{-\alpha}$
in a double power series as follows.

\begin{nr}\label{nr:stimaperron1}
For each $r\in (0,1]$ and $\alpha=2\beta > 0 $
\[
|1-r\xi|^{-\alpha} = \sum_{n=-\infty}^{+\infty} b_n \xi^n,
\]
with, for each $n\geq 0$,
\[
b_n  = b_{-n} = \dfrac{\sin( \beta \pi)}{ \pi }   r^n \int_0^1
(1-t)^{-\beta} t^{\beta-1}  t^n {(1-tr^2)}^{-\beta} \,dt.
\]
\end{nr}
\begin{proof}
\[
\begin{aligned}
|1-r\xi|^{-\alpha} &
= (1- r\xi)^{-\beta}(1-r\xi^{-1})^{-\beta} \\
&= \left( \sum_{k=0}^\infty \binom{-\beta}{k} (-r\xi)^k \right)
\cdot
\left( \sum_{h=0}^\infty \binom{-\beta}{h} (-r\xi^{-1})^h \right)\\
&= \sum_{h,k=0}^\infty \binom{-\beta}{k}  \binom{-\beta}{h}
(-r)^{k+h} \xi^{k-h}  \\ &= \sum_{n=-\infty}^\infty
\underbrace{\left( (-1)^n \sum_{\substack{k-h=n\\k,h\geq 0}}
\binom{-\beta}{k} \binom{-\beta}{h}  r^{k+h} \right)}_{ b_n } \xi^n
\end{aligned}.
\]
Now, recall that  for each $\beta>0$ and $N$ integer
\[\begin{aligned}
\binom{-\beta}{N} =  (-1)^N \binom{N+\beta - 1}{N} & = (-1)^N
\dfrac{\Gamma(N+\beta) \Gamma(1-\beta)}%
{\Gamma(N+1)\Gamma(\beta)\Gamma(1-\beta)}  \\
&= \dfrac{(-1)^N }{ \Gamma(\beta)\Gamma(1-\beta) } \cdot
\dfrac{\Gamma(N+\beta) \Gamma(1-\beta)}{\Gamma(N+1)} \\
&= \dfrac{(-1)^N \sin( \beta \pi)}{ \pi } \cdot
B ( 1-\beta, N+\beta ) \\
\end{aligned}\]
where $B(x,y)$ denotes the beta function, defined as
\[
B(x,y) = \int_0^1 t^{x-1} (1-t)^{y-1}\,dt =
\dfrac{\Gamma(x)\Gamma(y)}{\Gamma(x+y)}
\]
and we have used the equalities
\[
\Gamma(\beta)\Gamma(1-\beta) = \dfrac{\pi}{\sin (\beta \pi ) },
\]
\[
\binom{-\beta}{N} = \dfrac{(-\beta)(-\beta-1) \ldots (-\beta-N+1)}{
N!} = (-1)^N \dfrac{\Gamma(N+\beta)}{\Gamma(N+1)\Gamma(\beta)}
\]
and
\[
\binom{\beta}{N} =
\dfrac{\Gamma(\beta+1)}{\Gamma(N+1)\Gamma(\beta-N+1)}.
\]
We can now use the integral representation of the binomial function
\begin{equation}\label{eq:cinque}
\binom{-\beta}{N} = (-1)^N\dfrac{\sin( \beta \pi)}{ \pi } \int_0^1
(1-t)^{-\beta} t^{\beta-1} t^N \,dt,
\end{equation}
which implies that, by setting $N=h+n$,
\[\begin{aligned}
b_n  &= (-1)^n \sum_{h=0}^\infty \left( (-1)^{n+h}\dfrac{\sin( \beta
\pi)}{ \pi } \int_0^1 (1-t)^{-\beta} t^{\beta-1} t^{n+h} \,dt
 \binom{-\beta}{h}  r^{n+2h}
\right) \\
&= \dfrac{\sin( \beta \pi)}{ \pi }   r^n \sum_{h=0}^\infty \left(
(-1)^{h} \int_0^1 (1-t)^{-\beta} t^{\beta-1} t^{n+h} \,dt
 \binom{-\beta}{h}  r^{2h}
\right) \\
&= \dfrac{\sin( \beta \pi)}{ \pi }   r^n \int_0^1 (1-t)^{-\beta}
t^{\beta-1}  t^n \sum_{h=0}^\infty
\underbrace{%
\left( (-1)^{h} t^{h}
 \binom{-\beta}{h}  r^{2h} \right)}_{\binom{-\beta}{h} (-tr^2)^h}
\,dt
 \\
&= \dfrac{\sin( \beta \pi)}{ \pi }   r^n \int_0^1 (1-t)^{-\beta}
t^{\beta-1}  t^n {(1-tr^2)}^{-\beta} \,dt.
 \\
\end{aligned}
\]
\end{proof}

\begin{nr}
\label{nr:stimaperron2} For each $\beta\in (0,1)$, $r\in (0,1]$ and
integer $l\geq 2$
\[
\perron_l \left ( |1-r\xi|^{-\alpha} \right) =
\sum_{n=-\infty}^{+\infty} b_{ln} \xi^n.
\]
\end{nr}
\begin{proof}
It follows directly from equation \ref{eq:proprperron}. The
convergence is easy to check.
\end{proof}

\begin{nr}\label{nr:stimaperron3}
For each $\beta\in(0,1)$ and $r\in (0,1]$ and integer $l\geq 2$
\[
\perron_l \left ( |1-r\xi|^{-\alpha} \right) = \dfrac{\sin( \beta
\pi)}{ \pi } \int_0^1
\dfrac{(1-t)^{-\beta} t^{\beta-1} }{ %
{(1-tr^2)}^{\beta} } \dfrac{1 - (tr)^{2l}}{|1-(tr)^l\xi|^2} \,dt.
\]
\end{nr}
\begin{proof}
\[\begin{aligned}
\perron_l\left(|1-r\zeta|^{-\alpha} \right) \!\!\!&=
\sum_{n=-\infty}^\infty b_{ln} \zeta^n = \sum_{n=0}^\infty b_{ln}
\zeta^n +
\sum_{n=1}^\infty b_{ln} \zeta^{-n} \\
\!\!\!\!&=\!\!\! \sum_{n=0}^\infty \left( \dfrac{\sin( \beta \pi)}{
\pi } r^{ln} \int_0^1 (1-t)^{-\beta} t^{\beta-1}  t^{ln}
{(1-tr^2)}^{-\beta} \,dt \right) \zeta^n\!+
\\
\!\!\!& +\!\!\! \sum_{n=1}^\infty \left( \dfrac{\sin( \beta \pi)}{
\pi } r^{ln} \int_0^1 (1-t)^{-\beta} t^{\beta-1}  t^{ln}
{(1-tr^2)}^{-\beta} \,dt
\right) \zeta^{-n} \\
\!\!\!&=\!\!\! \dfrac{\sin( \beta \pi)}{ \pi }\!\! \int_0^1\!\!
\dfrac{(1-t)^{-\beta} t^{\beta-1} }{ %
{(1-tr^2)}^{\beta} } \left[ \sum_{n=0}^\infty (tr)^{ln} \zeta^n \!\!
+\!\! \sum_{n=1}^\infty (tr)^{ln} \zeta^{-n} \right]\!\!dt.
\end{aligned}\]
% But
% \[\begin{aligned}
% \sum_{n=0}^\infty
% (tr)^{ln} \zeta^n   &= \dfrac{1}{1-(tr)^l\zeta } \\
% \sum_{n=1}^\infty
% (tr)^{ln} \zeta^{-n} &=
% \dfrac{(tr)^l\zeta^{-1}}{1-(tr)^l\zeta^{-1}}, \\
% \end{aligned}
% \]
The conclusion follows since
\[
\sum_{n=0}^\infty (tr)^{ln} \zeta^n   + \sum_{n=1}^\infty (tr)^{ln}
\zeta^{-n}   = \dfrac{1 - (tr)^{2l}}{|1-(tr)^l\zeta|^2}.
\]
\end{proof}
Thus we proved the following result.
\begin{propo}
For $\beta\in(0,1)$, $r\in (0,1]$ and $\xi\in S^1\subset \CC$ the
potential $U$ can be written as
\[
U = |z|^{-\alpha} \left[ c_l  + l r^{-\beta} \dfrac{\sin( \beta
\pi)}{ \pi } \int_0^1
\dfrac{(1-t)^{-\beta} t^{\beta-1} }{ %
{(1-tr^2)}^{\beta} } \dfrac{1 - (tr)^{2l}}{|1-(tr)^l\xi^l|^2} \,dt
\right],
\]
where $c_l$ is the constant $c_l  = \sum_{j=1}^{l-1} | 1 - \zeta_l^j
|^{-\alpha}$.
\end{propo}

\begin{remark}\label{remark:qqq}
An analogue of the integral representation of the potential is well
known, and can be traced back to F.-F. Tisserand's book
\cite{Tisserand}  (chapter XVII)  for the exponent $\alpha=1$; it
had been used by M. Lindow \cite{Lindow24} (section 3) in computing
central configurations for the planar gravitational field  generated
by a regular $n$-gon. More recently D. Bang and B. Elmabsout
extended and generalised Lindow's theorem, proving an equivalent of
\ref{propo:potentialintegral} (Proposition 7 and 8 of
\cite{BanElm}). The proof given here is direct, and allows explicit
estimates that can be used to compute the Hessian for the potential
restricted to the shape sphere. Furthermore, it involves an
interesting connection with the $l$-adic Ruelle--Perron--Frobenius
operator (see P. Gaspard's paper \cite{Gaspard92}), which is worth a
mention.
\end{remark}

%%%=========================================================================

\end{document}